\documentclass[12pt]{article}
\usepackage[english]{babel}
\usepackage{amsmath,cite,amsthm}
\usepackage{amssymb,euscript}
\usepackage{latexsym}
\newtheoremstyle{theorem}%name
  {10pt}          % space above
  {10pt}  % space below
  {\sl}  % bofy font
  {\parindent}     % ident - empty=no indent,  \parindent= paragraph indent
  {\bf}  % thm head font
  {. }    % punctuation after thm head
  { }    % space after thm head: `` ``=normal \newline=linebreak
  {}     % thm head specification
\theoremstyle{theorem}
\newtheorem{theorem}{Theorem}

\newtheoremstyle{remark}%name
  {10pt}          % space above
  {10pt}  % space below
  {\sl}  % bofy font
  {\parindent}     % ident - empty=no indent,  \parindent= paragraph indent
  {\bf}  % thm head font
  {. }    % punctuation after thm head
  { }    % space after thm head: `` ``=normal \newline=linebreak
  {}     % thm head specification
\theoremstyle{remark}
\newtheorem{remark}{Remark}

\newtheoremstyle{defi}%name
  {10pt}          % space above
  {10pt}  % space below
  {\rm}  % bofy font
  {\parindent}     % ident - empty=no indent,  \parindent= paragraph indent
  {\bf}  % thm head font
  {. }    % punctuation after thm head
  { }    % space after thm head: `` ``=normal \newline=linebreak
  {}     % thm head specification
\theoremstyle{defi}

\newtheoremstyle{lemma}%name
  {10pt}          % space above
  {10pt}  % space below
  {\sl}  % bofy font
  {\parindent}     % ident - empty=no indent,  \parindent= paragraph indent
  {\bf}  % thm head font
  {. }    % punctuation after thm head
  { }    % space after thm head: `` ``=normal \newline=linebreak
  {}     % thm head specification
\theoremstyle{lemma}
\newtheorem{lemma}{Lemma}

\begin{document}

\title{New biharmonic bases\\ in commutative algebras of the second rank and monogenic functions\\ related  to the biharmonic equation}

\author{S.V. Gryshchuk\\
 Institute of Mathematics, \\ National Academy of Sciences
of Ukraine,\\ Tereshchenkivska Str. 3, 01024, Kyiv, Ukraine\\
gryshchuk@imath.kiev.ua, serhii.gryshchuk@gmail.com}

\maketitle

\begin{abstract}
Among all two-dimensional commutative algebras of
 the second rank  %over the field of complex numbers
  a totally of all their  biharmonic bases $\{e_1,e_2\}$, satisfying conditions
$\left(e_1^2+ e_2^2\right)^{2} = 0$, $e_1^2 + e_2^2 \ne 0$, is found in an explicit form.
A set of "analytic" (monogenic) functions  satisfying the biharmonic equation and defined in the real planes  generated by the biharmonic bases is built.  A characterization of biharmonic functions in bounded simply connected domains by real components  of some monogenic functions is found.

 \vspace*{6mm}
 {\bf AMS Subject Classification (2010):} Primary 30G35, 31A30; Secondary
 74B05 %2010

 % 74B05 Classical linear elasticity,  31A30 Biharmonic, polyharmonic functions and equations, Poisson's equation
%30G35 Functions of hypercomplex variables and generalized variables

  {\bf Key Words and Phrases:} biharmonic equation, biharmonic algebra, biharmonic base,  monogenic function of the biharmonic plane, isotropic plane strain.
\end{abstract}

\newpage
%%%%%%%%%%%%%%%%%%%%%%%%%%%%%%%%%%%%%%%%%%%%SECTION-1%%%%%%%%%%%%%%%%%%%%%%%%%%%%%%%%%%%%%
\section{Sta\-te\-ment of the Prob\-lem}                %{INTRODUCTION}

We say that an associative  algebra of the second rank
 with identity commutative over the field of complex (or real) numbers
$\mathbb C$ (\,$\mathbb{R}$\,)  is {\it biharmonic} if it contains  at least one
{\it biharmonic} basis, i.e., a basis $\{e_1,e_2\}$ satisfying the
conditions
\begin{equation}\label{biharm-bas}
 e_1^4 + 2e_1^2 e_2^2 + e_2^4 = 0,\qquad e_1^2+e_2^2\ne 0\,.
\end{equation}
%\EuScript{L}(e_1,e_2):=

A {\bf problem} of describing all biharmonic algebras and all their biharmonic  bases  has been posed by
I.\,P. Mel'nichenko in \cite{Mel86}.

In \cite{Mel86}, I.~P.~Mel'nichenko proved the uniqueness  of the biharmonic algebra
\begin{equation}\label{biharAlg}
  \mathbb{B}:=\{z_1\, e + z_2\, \rho:\, z_{k}\in \mathbb{C}, k=1,2\},\,  \rho^2 =0, e\rho =\rho, e^2 = e,
\end{equation}
here $e$ is the identity of algebra.
Note that the algebra \eqref{biharAlg} is isomorphic to the four-dimensional algebras over the field of real numbers considered by L.~Sobrero (cf., e.g., \cite{Sodbero}) and A.~Douglis \cite{Douglis-53}.

V.\,F. Kovalev and I.\,P. Mel'nichenko \cite{Kov-Mel}
 found a
multiplication table for a biharmonic basis $\{e_1,e_2\}$:
\begin{equation} \label{tab_umn_bb}
e_1=e,\qquad e_2^2=e_1+2ie_2,
\end{equation}
where $i$\,\, is the imaginary complex unit.

I.~P.~Mel'nichenko claimed that he found a complete totally of biharmonic bases and described this  set of biharmonic bases in  \cite{Mel86} also.

There were considered in
\cite{GrPl_umz-09,Cont-13} %KM-BFf,UMZ_mon_09}
 some $\mathbb{B}$-valued "analytic" functions $\Phi$, defined in domains which lies in the span of the vectors  \eqref{tab_umn_bb} over the field of real numbers, %$\mathbb{R}$: $\{xe_1 + ye_2\}$
such that real components of function  $\Phi$ (over the system of vectors: $e_k$, $ie_k$, $k=1,2$) satisfy the biharmonic equation
\begin{equation} \label{gen-bih-eq}
\left(\Delta_2\right)^{2} u(x,y):= \left(\frac{\partial^4}{\partial x^4} + 2 \frac{\partial^4
}{\partial x^2\partial y^2}+
\frac{\partial^4 }{\partial y^4} \right) u(x,y)=0, %\qquad e_1^2+e_2^2\ne 0,
\end{equation}
where
$\Delta_2:=\frac{\partial^2}{\partial x^2}+ \frac{\partial^2}{\partial y^2}$, i.e., these components are biharmonic functions  (in some domain of the Cartesian plane  $xOy$ which depends on the domain where a function $\Phi$ is defined).

It is proved in  \cite{GrPl_umz-09}
 that for  any biharmonic function  $u$, defined in the
  bounded and simply connected domain,
  exist functions  $\Phi$ such that their fixed component (before  $e_1$) coincides with  $u$, and, described the whole totally of these functions  $\Phi$.

  Consider a short summary of new results which are described  in the proposed paper.
   A totally of biharmonic bases described in  \cite{Mel86} is not complete,  and the whole totaly is found also here. A constructive description of monogenic functions, defined in domains belong to the planes generated by any biharmonic bases, is found in an explicit form. It is proved
 that for  any biharmonic function  $u$, defined in the
  bounded and simply connected domain,
  exist some monogenic  functions  $\Phi$ such that their fixed component (before  $e_1$) coincides with  $u$, and, described the whole totally of these functions  $\Phi$.

  %%%%%%%%%%%%%%%%%%%%%%%%%%%%%%%%%%%%%%%%%%%%%%%%%%%%%%%%%%%%%%%%%SECTION-2%%%%%%%%%%%%
  \section{Biharmonic bases} Here and what follows we assume  to choose  an upper or lower sign in formulas containing the symbol $\pm$.
The next teorem gives a constructive description of all biharmonic bases.

\begin{theorem}\label{all-basesBiharm}
 \emph{All biharmonic bases (of the algebra
$\mathbb{B}$) %, ùî çàäîâîëüíÿþòü óìîâó \emph{(\ref{tab_umn_ba-1})}
  are expressed by the formulas\emph{}:
\begin{equation}\label{E-12B-p1}
  e_1=\alpha_1\, e +\beta_1\, \rho,\,
e_{2}=\pm i \left(\alpha_1\,e + \beta_2\, \rho \right),
\end{equation}
where complex numbers  $\alpha_k$, $\beta_k$, $k=1,2$,  %îäíî÷àñíî
satisfy conditions\emph{:} %÷îòèðüîõ
$\alpha_1 \ne 0$, $\beta_1 \ne \beta_2$.
 }\end{theorem}

This theorem can be proved by similar arguments which were done at the proof of Theorem~1 in \cite{Gr-umzOrth18-1-An}.

\begin{remark}
Biharmonic bases of the form
\eqref{E-12B-p1} with $\beta_2 = \beta_1 -\frac{1}{2\alpha_1}$ \emph{(}this and only this subset of all biharmonic bases \eqref{E-12B-p1}\emph{)} are considered in
\emph{\cite{GrPl_umz-09}}.
\end{remark}

%%%%%%%%%%%%%%%%%%%%%%%%%%%%%%%%%%%%%%%%%SECTION-3%%%%%%%%%%%%%%
\section{Monogenic functions of the plane generated by elements of the biharmonic bases}

Consider a plane  $\mu_{e_1, e_2}:=\{xe_1+ye_2:
x,y\in\mathbb{R}\}$ where  $\{e_1,e_2\} $ is a fixed
bases in \eqref{E-12B-p1}.
We use  the Euclidean norm $\|a\|:=\sqrt{|z_1|^2+|z_2|^2}$, where
$a=z_1e_1+z_2e_2\in\mathbb{B}$, $z_k \in \mathbb{C}$, $k=1, 2$. With a domain  $D$
of the Cartesian plane  $xOy$  we associate the congruent domain
 $D_{\zeta}:=
\{\zeta = x e_1+y e_2~\in~\mu_{e_1, e_2}: (x,y)\in D\} \subset \mu_{e_1,e_2}$.

In what follows,
\[(x,y)\in D, \, \zeta=x e_1 + y e_2 \in D_{\zeta}, \, z=x+iy\in \mathbb{C}\equiv \mathrm{Re\,}z +
 i\mathrm{Im\,}z; \, x, y \in \mathbb{R}.\]

Note that every element $\zeta\in\mu_{e_1, e_2}$, $\zeta\ne 0$, is invertible  (exists an element $\zeta^{-1}\in\mathbb{B}$: $\zeta \zeta^{-1}=e$).

Inasmuch as divisors of zero don't belong to the plane $\mu_{e_1,e_2}$,
one can define the derivative $\Phi'(\zeta)$ of function $\Phi
\colon D_{\zeta}\longrightarrow \mathbb{B}$ in the same way as in
the complex plane:
\[\Phi'(\zeta):=\lim\limits_{h\to 0,\, h\in\mu_{e_1, e_2}}
\bigl(\Phi(\zeta+h)-\Phi(\zeta)\bigr)\,h^{-1}\,.\]
 We say that a function $\Phi : D_{\zeta}\longrightarrow \mathbb{B}$ is
\textit{monogenic} in a domain $D_{\zeta}$  iff the derivative
$\Phi'(\zeta)$ exists at  every point $\zeta\in D_{\zeta}$.

%%%%%%%%%%%%%%%%%%%%%%%%%%%%%%%%%%%%%
%%%%%%%%%%%%%%%%%%%%%%%%%%%%%%%%%%%%%
 %We use the notion of monogenic function in the sense of
%existence of derived numbers for this function in the domain
%$D_{\zeta}$ (cf., e.g., \cite{Goursat,Trokhimchuk}).

Every function  $\Phi \colon D_{\zeta}\longrightarrow \mathbb{B}$
has a form
\begin{equation}\label{mon-funk} \Phi(\zeta)=U_{1}(x,y)\,e_1+U_{2}(x,y)\,ie_1+
U_{3}(x,y)\,e_2+U_{4}(x,y)\,ie_2,
\end{equation}
where $\zeta=xe_1+ye_2$, $U_{k}\colon D\longrightarrow \mathbb{R}$,
$k=\overline{1,4}$.

Every real component  $U_{k}$, $k=\overline{1,4}$, in expansion
 (\ref{mon-funk}) we denote by
$\mathrm{U}_{k}\left[\Phi\right]$, i.e., for %arbitrary fixed
$k\in\{1,\dots,4\}$:
$\mathrm{U}_{k}\left[\Phi(\zeta)\right]:=U_{k}(x,y)$ for all $\zeta=
xe_1 + y e_2 \in D_{\zeta}$.

It follows from \eqref{E-12B-p1} formulas:
\begin{equation}\label{EkBiktTablUmn}
  (e_{k})^{2} = (-1)^{k+1}  \alpha_{1} \left(\alpha_1\, e + 2 \beta_k \, \rho \right), \, k=1,2,\,
e_1 e_2 = \pm i \alpha_1 \left(\alpha_1 e + \left(\beta_1 + \beta_2\right)\right)  \rho.
\end{equation}
 In \eqref{EkBiktTablUmn} and till the end of the paper by the double sing  $\pm$  we shall mean the same single  sign as in  \eqref{E-12B-p1} for corresponding basis  $\{e_1, e_2\}$.

 Analogously to the case of the  biharmonic basis \eqref{tab_umn_bb} (see \cite{Kov-Mel,GrPl_umz-09}) we obtain
the next theorem. %by use of\eqref{EkBiktTablUmn}

\begin{theorem}\label{CRim-BGenBikh}
 A function
$\Phi \colon  D_{\zeta}\longrightarrow \mathbb{B}$ is monogenic in a
domain $D_{\zeta}$ if and only if components $U_{k}$,
$k=\overline{1,4}$, of the expression \eqref{mon-funk} are
differentiable in the domain $D$ and the following analog of the
Cauchy -- Riemann conditions  is fulfilled\emph{:}
\begin{equation}\label{usl_K_R}
\frac{\partial \Phi(\zeta)}{\partial y}\,e_1=\frac{\partial
\Phi(\zeta)}{\partial x}\,e_2.
\end{equation}
\end{theorem}

\begin{remark}
Using \eqref{EkBiktTablUmn}, we obtan that
in an extended form the condition \eqref{usl_K_R} for the monogenic
function \eqref{mon-funk} is equivalent to the system of four
equations  with respect to
components $U_{k}=\mathrm{U}_{k}\left[\Phi\right]$,
$k=\overline{1,4}$, in \eqref{mon-funk}:
\begin{eqnarray}
\frac{\partial U_{1}(x,y)}{\partial y}- \left(\pm\right)\frac{\partial U_{4}(x,y)}{\partial y} \pm \frac{\partial U_{2}(x,y)}{\partial x}+
 \frac{\partial U_{3}(x,y)}{\partial x} = 0, \label{kr1-0BikhB-1}\\
\frac{\partial U_{2}(x,y)}{\partial y} \pm \frac{\partial U_{3}(x,y)}{\partial y} - \left(\pm \right) \frac{\partial U_{1}(x,y)}{\partial x}+
 \frac{\partial U_{4}(x,y)}{\partial x} = 0, \label{kr1-0BikhB-2}\\
2 \mathrm{Re}\beta_{1}\,\frac{\partial U_{1}(x,y)}{\partial y} - 2 \mathrm{Im}\beta_{1}\,\frac{\partial U_{2}(x,y)}{\partial y} -
 \left(\pm \mathrm{Im}\left(\beta_1 + \beta_2\right)\right)
  \frac{\partial U_{3}(x,y)}{\partial y}-\nonumber \\
-\left(\pm \mathrm{Re} \left(\beta_1 + \beta_2\right)\right)
\frac{\partial U_{4}(x,y)}{\partial y}
 \pm \mathrm{Im} \left(\beta_1 + \beta_2\right) \frac{\partial U_{1}(x,y)}{\partial x}+
 \nonumber\\
 \pm \mathrm{Re} \left(\beta_1 + \beta_2\right)
  \frac{\partial U_{2}(x,y)}{\partial x}
+ 2 \mathrm{Re}{\beta_2} \frac{\partial U_{3}(x,y)}{\partial x} -
2 \mathrm{Im}{\beta_2} \frac{\partial U_{4}(x,y)}{\partial x}&=&0, \label{kr1-0BikhB-3} \\
2 \mathrm{Im}\beta_{1}\,\frac{\partial U_{1}(x,y)}{\partial y} +
 2 \mathrm{Re}\beta_{1}\,\frac{\partial U_{2}(x,y)}{\partial y}
  \pm \mathrm{Re}\left(\beta_1 + \beta_2\right) \frac{\partial U_{3}(x,y)}{\partial y}-\nonumber \\
-\left(\pm\mathrm{Im} \left(\beta_1 + \beta_2\right)\right)
\frac{\partial U_{4}(x,y)}{\partial y} -
\left(\pm \mathrm{Re} \left(\beta_1 + \beta_2\right) \right)\frac{\partial U_{1}(x,y)}{\partial x}+
 \nonumber\\
 \pm  \mathrm{Im} \left(\beta_1 + \beta_2\right)
\frac{\partial U_{2}(x,y)}{\partial x}+
 2 \mathrm{Im}{\beta_2} \frac{\partial U_{3}(x,y)}{\partial x} +
2 \mathrm{Re}{\beta_2} \frac{\partial U_{4}(x,y)}{\partial x}&=&0. \label{kr1-0BikhB-4}
\end{eqnarray}
\end{remark}

Take into consideration a variable and a domain of the complex plane of the form:
\begin{equation}\label{DBikharm}
Z= \alpha_1 \left(x\pm iy\right),\,
  D_{Z}:=\{Z=\alpha_1 \left(x\pm iy\right): (x,y)\in D\}.
\end{equation}

Similar to the proof of theorems 1 and 2 in \cite{GrPl_umz-09} with use of Theorem~\ref{CRim-BGenBikh}, we obtain  an expression of monogenic functions via holomorphic functions of complex variable $Z$ in the domain $D_{Z}$.

\medskip
 \begin{theorem}\label{gCR-holom-Biharm}
  A function
 $\Phi \colon D_{\zeta} \longrightarrow \mathbb{B}$ is monogenic in  $D_{\zeta}$
 if and only if the following equality if fulfilled
\begin{equation}\label{rep-mon-f-B-0Biharm}
\Phi(\zeta)=
F(Z)\, e+
\left(\left(\frac{\beta_1}{\alpha_1}\, Z \pm i\left(\beta_2 - \beta_1\right)y\right)\, F'(Z) + F_{0}(Z)\right)\rho
\equiv  \Phi\left[F, F_{0}\right](\zeta)
 \,\, \forall \zeta \in D_{\zeta},
\end{equation}
where $F$, $F_{0}$ are some holomorphic functions of the complex variable
$Z$ in the domain  $D_{Z}$, $F'$ if the derivative of  $F$.
\end{theorem}

In what follows, we mean by $\Phi\left[F, F_{0}\right]$ an expression of the form \eqref{rep-mon-f-B-0Biharm}.

%\medskip
\begin{remark}
An expression  \eqref{rep-mon-f-B-0Biharm} for monogenic functions with
biharmonic basis
\eqref{tab_umn_bb} is found in  \emph{\cite{GrPl_umz-09}}, and before, for some particular class of domains  $D_{\zeta}$~--- in the paper \emph{\cite{Kovalov}}.
\end{remark}

%\medskip
\begin{remark}
A relation between monogenic functions   \eqref{tab_umn_bb}
\emph{(}related to the biharmonic basis
\eqref{tab_umn_bb}\emph{)} with analytic by Douglis functions \emph{(}considered, for example, in \emph{\cite{Sold-hyperan-SMiPr-04,Sold-SMFN-dr_nar-16}}\emph{)}
is found in \emph{\cite{Gr-umzOrth18-1-An}}.
\end{remark}

%\medskip
\begin{remark}
Consider the biharmonic base $\{e_1, e_2\}$ such that  $e_2=e$. Then in  \eqref{E-12B-p1} we have
\begin{equation}\label{bikhE}
  e_1 = e,\,\, e_2 =(\pm i) e + (\pm i\beta_2)\rho, \,
   \beta_2 \in \mathbb{C}, \beta_2 \ne 0.
\end{equation}
  Thus,
\eqref{gCR-holom-Biharm} turns onto
\begin{equation}\label{rep-mon-f-B-0Bih-E}
\Phi(\zeta)=
F(x\pm iy)\, e+
\left(\pm i \beta_2 y \, F'(x\pm iy) + F_{0}(x\pm iy)\right)\rho
 \,\, \forall \zeta \in D_{\zeta}.
 \end{equation}

The expression  \eqref{rep-mon-f-B-0Bih-E} is obtained for  ``monogenic functions'' $\Phi \colon D_{\zeta} \longrightarrow \mathbb{B}$, which are understood  as continuous functions and differentiable by Gateaux (along the positive rays) in the paper \emph{\cite{Sh-zb-15}} \emph{(}a fact that a  ``convexity''  of the domain $D_{\zeta}$ is omitted can be proved analogous to the similar fact in \emph{\cite[after Theorem~3]{{Gr-umzOrth18-1-An}}}\emph{)}.
 Taking into account  that monogenic function  \eqref{rep-mon-f-B-0Bih-E} is continuous also, obtain that these two kinds of monogeneity  are equivalent in the case of the biharmonic bases  \eqref{bikhE}.
\end{remark}
%%%%%%%%%%%%%%%%%%%%%%%%%%%%%%%%%%%%%%%%%%%%%

\vspace*{2mm} {\bf 4.  Biharmonic functions as components of monogenic functions.}
 The equality \eqref{rep-mon-f-B-0Biharm} and Theorem~\ref{CRim-BGenBikh} yield that any monogenic function  $\Phi \colon D_{\zeta} \longrightarrow \mathbb{B}$  has a derivative  $\Phi^{(n)}$ for every natural $n$, and
  components  $U_{k}=\mathrm{U}_{k}\left[\Phi\right]$, $k=\overline{1,4}$, are infinitely differentiable in $D$. Therefore, using the equality
\begin{equation}\label{bikhMon}
  \left(\Delta_{2}^{2}\right)\Phi(\zeta)= \EuScript{L}(e_1,e_2)\Phi^{(4)}(\zeta)\equiv 0 \,\, \forall \zeta \in D_{\zeta},
\end{equation}
we conclude that components  $U_{k}$, $k\in\{1,\dots,4\}$, from the expansion  \eqref{mon-funk} are
biharmonic functions in the domain  $D$.

 Now and till the end of the paper we mean by the biharmonic basis \eqref{E-12B-p1} the following basis:
\begin{equation}\label{BihhOneNew}
  e_1 = e + \rho,\,  e_2 = i\left(e+ 2 \rho\right).
\end{equation}
Note, that the basis   \eqref{BihhOneNew} was not found in  \cite{Mel86}.
Then in \eqref{DBikharm} we take $\alpha=1$ and the sigh ``+''.
 Thus,  the formula \eqref{rep-mon-f-B-0Biharm}  turns onto the form
\begin{equation}\label{rep-monf-BiharmNewB}
\Phi(\zeta)=
F(z)\, e+
\left(\left(z+iy\right)\, F'(z) + F_{0}(z)\right)\rho
 \,\, \forall \zeta \in D_{\zeta}.
\end{equation}

 With respect to the basis
\eqref{BihhOneNew} the expression
\eqref{rep-mon-f-B-0Biharm}
turns onto the form
\begin{equation}\label{rep-mon-f-B-0BiharmE}
  \Phi(\zeta) =\left(F(z) +\widetilde{F}(z)\right)e_1 + \widetilde{F}(z)\,ie_2\,\, \forall \zeta \in D_{\zeta},
\end{equation}
here
$\widetilde{F}(z): = F(z) - \left(z+iy\right)F'(z)-F_{0}(z)$ for every  $z \in D_{z}$.

\medskip
\begin{lemma}\label{oneCompZeroBikhN}
Let the domain  $D$ be bounded and simply connected.
Any monogenic function  $\Phi_{0} \colon D_{\zeta} \longrightarrow \mathbb{B}$ such that
\begin{equation}\label{formUnZe}
  \mathrm{U}_{1}\left[\Phi_{0} (\zeta)\right] =0 \,\, \forall \zeta \in D_{\zeta},
\end{equation}
is expressed in the form\emph{:}
\begin{equation}\label{PhiZero}
  \Phi_{0}(\zeta)= a\left(\zeta +z\rho\right)  +  bie_1 + c e_2 + d ie_2,
\end{equation}
where  $a, b, c, d$ are any real constants.
\end{lemma}

\medskip
\begin{remark}
 An analogous result for
monogenic functions generated by the biharmonic basis \eqref{tab_umn_bb} is considered in \emph{\cite[Lemma~3]{GrPl_umz-09}}. Note, that in this case real components of an analog of the formula \eqref{PhiZero} are polynomials of degree no more then the second power with respect to real variables $x$ and  $y$,
unlike the formula \eqref{PhiZero} where components are  polynomials of degree no more then the first power with respect to real variables $x$ i $y$.
\end{remark}

\medskip
\begin{lemma}\label{vspomLemBikhNew}
 Let $f_{k} \colon D_{z} \longrightarrow \mathbb{C}$, $k=1,2$, are holomorphic functions of the complex variable. Then a monogenic function  $\Phi\colon D_{\zeta} \longrightarrow \mathbb{B}$ defined by the formula
\begin{equation}\label{vspomMonBikh}
  \Phi(\zeta):=\Phi\left[F, F_{0}\right](\zeta)\, \, \forall \zeta \in D_{\zeta},
\end{equation}
where %àíàë³òè÷í³ ôóíêö³¿ $F_{k}, F_{k0} \colon D_{z} \longrightarrow \mathbb{C}$
\begin{equation}
\begin{aligned}
F(z) := 2 f_{2}(z),
 F_{0}(z):= 4f_{2}(z)-f_{1}(z)-3 z f'(z)\,\, \forall z \in D_{z},
 \end{aligned}
 \end{equation}
satisfies the relation:
 \begin{equation}\label{Phiuno-tre}
   \mathrm{U}_{1}\left[\Phi(\zeta)\right]=\mathrm{Re}\, \left(f_{1}(z)+\overline{z}f_{2}'(z)\right)\, ( \overline{z}:= x-iy) \,\, \forall \zeta \in D_{\zeta}.
 \end{equation}
 \end{lemma}

 A proof of Lemma is a simple corollary of expressions  \eqref{rep-mon-f-B-0Biharm} and \eqref{rep-mon-f-B-0BiharmE}.

The well-known fact is that any biharmonic function  $U_{1}$ in %boundedâ îáìåæåí³é
 simply connected domain  $D$ is expressed by the Goursat formula  (cf., e.g., \cite{Mush_upr})
\begin{equation}\label{formGursa}
  U_{1}(x,y)=\mathrm{Re\,}\left(\psi(z) + \overline{z}\varphi(z)\right),
\end{equation}
where  $\psi$, $\varphi$ are holomorphic in  $D_{z}$  functions, $\overline{z}:= x-iy$,.

Combining the previous results we get the following theorem.

\medskip
 \begin{theorem}\label{prMonBikhN}
  Let a function  $U_{1}$ is biharmonic in the bounded and simply connected domain    $D$
   of the Cartesian plane  $xOy$. Then all monogenic functions $\Phi \colon D_{\zeta} \longrightarrow \mathbb{B}$ such that
 $\mathrm{U}_{1}\left[\Phi\right] \equiv U_{1}$ are expressed by the form
 \begin{equation}\label{prBikhMon}
   \Phi(\zeta)= \Phi\left[2\EuScript{F}, 4\EuScript{F} - \psi-3z \varphi\right](\zeta) + \Phi_{0}(\zeta) \,\, \forall \zeta \in D_{\zeta},
 \end{equation}
  where $\Phi_{0}$ is defined by the formula  \eqref{PhiZero},
  $\EuScript{F}$ is a primitive function for  $\varphi$.
  \end{theorem}

  \medskip
\begin{remark}
 An  analogous result for monogenic functions generated by the biharmonic basis  \eqref{tab_umn_bb}
 was found in  \emph{\cite[Theorem~5]{GrPl_umz-09}}.
\end{remark}

\medskip
\begin{remark}
An  analogous results for monogenic functions related to some bases generated by the stress equations
for certain orthotropic deformations  were found in
 \emph{\cite{Gr-umzOrth18-2-An,SlovZb18}}.
\end{remark}

\section{Acknowledgements}
The research is partially supported by the Fundamental Research Programme funded by the Ministry of Education  and Science  of Ukraine  (project No.~0116U001528) and
and a common grant of the National Academy of Sciences of Ukraine and the Poland Academy
of Sciences (No. 32/2018).

%"Topological-analytical methods in complex and hypercomplex analysis" of the Polish Academy of Sciences and the %National Academy of Sciences of Ukraine (2014 - 2017)

%\newpage

%%%%%%%%%%%%%%%%%%%%%%%%%%%%%%%%%%%%%%%%%%%%%%

%%%%%%%%%%%%%%%%%%%%%%%%%%%%%%%%%%%%%
%%%%%%%%%%%%%%%%%%%%%%%%%%%%%%%%%%%%%

%Results of this paper have been announced at a preprint of Arxiv...
%(arXiv: 1202.0993v1 [math.CV]).

%\section{Acknowledgements}

%The authors are grateful to Professor Massimo Lanza de Cristoforis
%for invitations and given good possibilities for their joint
%researches during the {\em{Mini-courses in Mathematical
%Analysis}}~--~2011.

%%%%%%%%%%%%


\newpage

\begin{thebibliography}{99}

\bibitem{Mel86}
 Mel'nichenko~I.\,P. {\em Biharmonic bases in algebras of the second rank},
{Ukr. Math. J.}, \textbf{38} (1986), No.~2, 252--254.

\bibitem{Kov-Mel}
 Kovalev~V.\,F.   Mel'nichenko~I.\,P. {\em Biharmonic functions on the
biharmonic plane}, {Reports Acad. Sci. USSR, ser. A.}, No.~8
(1981), 25--27 [in Russian].

 \bibitem{Sodbero}
 Sobrero~L. {\em   Nuovo metodo per lo studio dei problemi di
elasticit\`{a}, con applicazione al problema della piastra forata}.
Ricerche di Ingegneria. \textbf{13} (1934), No.~2, 255--264 [in Italian].

\bibitem{Douglis-53}
 Douglis~A., {\em  A function-theoretic approach to elliptic systems of equations in
two variables}, Communications on Pure and Applied Mathematics,
 \textbf{6} (1953), No.~2, 259--289.




%\bibitem{Kravchenko-Shapiro}
 %V.V. Kravchenko and  M.V. Shapiro, {\em Integral representations for
%spatial models of mathematical physics}, Pitman Research Notes in
%Mathematics, Addison Wesley Longman Inc. (1996).

%\bibitem{Spros}
% W. Spr\"{o}ssig, {\em Quaternionic analysis and Maxwell's
%equations}, CUBO A Math. J., {\bf 7}, no. 2 (2005), 57 -- 67.
%%%%%%%%%%%%%%%%%%%%%%%%%%%%%%%%%%%%%%%%%%%%%%%%%%%%%%%%%%%%%%%%%%%%%%%%%%%%%%%%%mine
\bibitem{GrPl_umz-09} Grishchuk~S.\,V. (the same as:  Gryshchuk~S.\,V.), Plaksa~S.\,A.
 Monogenic functions in a
biharmonic algebra, {Ukr. Mat. Zh.}, \textbf{61} (2009),
No.~12, 1587--1596 [in Russian]; {\bf English translation:} % (Springer) in
{Ukr. Math.~J.}, \textbf{61} (2009), No.~12, 1865--1876.

\bibitem{Cont-13} Gryshchuk~S.\,V., Plaksa S.\,A. {\em Basic Properties of Monogenic Functions
in a Biharmonic Plane}, in: \textit{``Complex Analysis and Dynamical
Systems V'',  Contemporary Mathematics}, \textbf{591} (2013),
 Amer. Math. Soc., Providence, RI, 127--134.
%(http://dx.doi.org/10.1090/conm/591/11831).

\bibitem{Gr-umzOrth18-1-An}
 Gryshchuk~S.\,V.  {\em   Commutative ñomplex algebras of the second rank with unity and some cases of plane orthotropy~I },  Ukr. Mat. Zh.,  \textbf{70} (2018), No.~8,  %¹ 8.
1058--1071 [in Ukrainian];
{\bf English translation:} {Ukr. Math. J.}, \textbf{70} (2019), No.~8, 1221--1236.
%Gryshchuk~S.\,V. \emph{Commutative Complex Algebras of the Second Rank with Unity and Some Cases of Plane %Orthotropy~I.}

\bibitem{Gr-umzOrth18-2-An}
 Gryshchuk~S.\,V.  {\em   Commutative ñomplex algebras of the second rank with unity and some cases of plane orthotropy~II },  Ukr. Mat. Zh.,  \textbf{70} (2018), No.~10,  %¹ 8.
1382--1389 [in Ukrainian]; {\bf English translation:} {Ukr. Math. J.}, \textbf{70} (2019), No.~10, 1594--1603.

\bibitem{Kovalov}
 V. F.~Kovalev,  \textit{Biharmonic Schwarz problem}, Preprint
No.~86.16.~ Institute of Mathematics, Acad. Sci. USSR, Inst. of
Math. Publ. House, Kiev, 1986 [in Russian].

\bibitem{Sold-hyperan-SMiPr-04}
     Soldatov~A.\,P. {\em Hyperanalytic functions and their applications}, J. of Math. Sci.
      \textbf{132} (2006), No.~6 , 827--882.


\bibitem{Sold-SMFN-dr_nar-16}  Soldatov~A.\,P. {\em On the theory of anisotropic plane elasticity},  Contemporary Mathematics. Fundamental Trends,  RUDN, Moscow,  \textbf{60} (2016), pp.~114--163 [in Russian].
     %Òðóäû Ñåäüìîé Ìåæäóíàðîäíîé êîíôåðåíöèè ïî äèôôåðåíöèàëüíûì è ôóíêöèîíàëüíî-äèôôåðåíöèàëüíûì óðàâíåíèÿì %(Ìîñêâà, 22–29 àâãóñòà, 2014). ×.~3.~---

 \bibitem{Sh-zb-15}Shpakivskyi~V.\,S.
{\em Monogenic functions in finite-dimensional commutative associative algebras},
Zb. Pr. Inst. Mat. NAN Ukr.,
{\bf 12} (2015), No.~3, 251--268.

\bibitem{Mush_upr}
 Muskhelishvili N.I. {\em Some basic problems of the mathematical theory of elasticity.
Fundamental equations, plane theory of elasticity, torsion and
bending. English transl. from the 4th Russian edition by
R.M.~Radok},  Noordhoff International Publishing: Leiden, 1977.

\bibitem{SlovZb18}
Gryshchuk~S.\,V. {\em Monogenic functions in two dimensional commutative algebras to equations of plane
orthotropy}, Pratsi Inst. Prikl. Mat. Mikh. NANU (Proc. of In-te of Appl. Math. and Mech. of NAS of Ukraine),  \textbf{32} (2018), %Slovyansk: "TexPrintCentre", pp.~
18--29 [in Ukrainian].
%%%%%%%%%%%%%%%%%%%%%%%%%%%%%%%%%%%%%%%%%%%%%%%%%%%%%%%%%%%%%%%%%%%%%%%%%%%%%%%%%%END of PUbl List


\end{thebibliography}
\end{document}